\documentclass[11pt,a4paper]{amsart} 
\usepackage[all]{xy}
\usepackage{ifthen} 
\usepackage{amssymb}
\calclayout
\makeatletter 
\def\serieslogo@{} 
\def\@setcopyright{} 
\makeatother

\title{Support varieties -- an ideal approach}

\author{Aslak Bakke Buan}
\address{Aslak Bakke Buan\\ 
Institutt for matematiske fag\\ NTNU\\
N-7491 Trondheim\\ Norway.} 
\email{aslakb@math.ntnu.no}
\author{Henning Krause}
\address{Henning Krause\\ Institut f\"ur Mathematik\\
Universit\"at Paderborn\\ 33095 Paderborn\\ Germany.}
\email{hkrause@math.uni-paderborn.de}
\author{{\O}yvind Solberg} 
\address{{\O}yvind Solberg\\
Institutt for matematiske fag\\ NTNU\\
N-7491 Trondheim\\ Norway.} 
\email{oyvinso@math.ntnu.no}

\newtheorem{lem}{Lemma}[section]
\newtheorem{prop}[lem]{Proposition} \newtheorem{cor}[lem]{Corollary}
\newtheorem{thm}[lem]{Theorem}

\theoremstyle{remark}

\newtheorem{rem}[lem]{Remark}
\newtheorem{exm}[lem]{Example}

\theoremstyle{definition}
\newtheorem{defn}[lem]{Definition}

\numberwithin{equation}{section}

\renewcommand{\mod}{\operatorname{mod}\nolimits}
\renewcommand{\leq}{\leqslant}
\renewcommand{\geq}{\geqslant}

\newcommand{\card}{\operatorname{card}\nolimits}

\newcommand{\fp}{\operatorname{fp}\nolimits}
\newcommand{\Mod}{\operatorname{Mod}\nolimits}
\newcommand{\End}{\operatorname{End}\nolimits}
\newcommand{\Hom}{\operatorname{Hom}\nolimits}

\newcommand{\supp}{\operatorname{supp}\nolimits}

\newcommand{\Ext}{\operatorname{Ext}\nolimits}
\newcommand{\HH}{\operatorname{HH}\nolimits}

\newcommand{\umod}{\operatorname{\underline{mod}}\nolimits}

\newcommand{\Inj}{\operatorname{Inj}\nolimits}
\newcommand{\Spec}{\operatorname{Spec}\nolimits}
\newcommand{\Proj}{\operatorname{Proj}\nolimits}
\newcommand{\Qcoh}{\operatorname{Qcoh}\nolimits}
\newcommand{\coh}{\operatorname{coh}\nolimits}

\newcommand{\Ab}{\mathrm{Ab}}

\newcommand{\op}{\mathrm{op}}
\newcommand{\inc}{\mathrm{inc}}

\newcommand{\can}{\mathrm{can}}

\newcommand{\id}{\mathrm{id}}

\newcommand{\per}{\mathrm{per}}
\newcommand{\thick}{\mathrm{thick}}
\newcommand{\open}{\mathrm{open}}

\newcommand{\comp}{\mathop{\raisebox{+.3ex}{\hbox{$\scriptstyle\circ$}}}}
\newcommand{\lto}{\longrightarrow}
\newcommand{\xto}{\xrightarrow}

\def\li{\varinjlim}

\def\a{\alpha}
\def\b{\beta}

\def\d{\delta}

\def\p{\phi}

\def\s{\sigma}
\def\t{\tau}

\def\La{\Lambda}
\def\Si{\Sigma}

\def\Om{\Omega}
\def\A{{\mathcal A}}
\def\B{{\mathcal B}}
\def\C{{\mathcal C}}
\def\D{{\mathcal D}}
\def\E{{\mathcal E}}

\def\L{{\mathcal L}}

\def\OO{{\mathcal O}}
\def\P{{\mathcal P}}

\def\Y{{\mathcal Y}}

\def\bbZ{\mathbb Z}

\def\bfD{\mathbf D}

\def\bfL{\mathbf L}

\def\fp{\mathfrak p}

\begin{document}

\begin{abstract}
We define support varieties in an axiomatic setting using the prime
spectrum of a lattice of ideals. A key observation is the
functoriality of the spectrum and that this functor admits an
adjoint. We assign to each ideal its support and can classify ideals
in terms of their support.  Applications arise from studying abelian
or triangulated tensor categories.  Specific examples from algebraic
geometry and modular representation theory are discussed, illustrating
the power of this approach which is inspired by recent work of Balmer.
\end{abstract}

\maketitle 
\setcounter{tocdepth}{1}
\tableofcontents

\section*{Introduction}

The spectrum of prime ideals and the support of objects like modules,
sheaves, complexes etc.\ belong to the fundamental concepts from
algebraic geometry. In fact, the use of these concepts is not
restricted to algebraic geometry and similar notions exist for
instance in modular representation theory.  In this paper we discuss a
general approach which allows us to study prime ideal spectra and
supports in various settings. 

A prime ideal spectrum comes naturally equipped with a topology which
is usually called Zariski topology.  However, there are various
instances where it is more natural to consider another `opposite'
topology. It is one of the principal aims of this work to clarify the
parallel use of two different topologies on a prime ideal
spectrum. This is based on the notion of a spectral space, first
introduced by Hochster \cite{Ho}.

We give a couple of motivating examples which illustrate the use of
such different topologies. Let $A$ be a commutative ring and denote by
$\Spec A$ the set of prime ideals, together with the usual Zariski
topology. We obtain another topology and write $\Spec^*A$ if we take
the quasi-compact Zariski open sets as a basis of closed sets. Now
suppose that $A$ is noetherian.  Let $\mod A$ denote the abelian
category of all finitely generated $A$-modules. Then the assignment
$$\mod A\supseteq \C\mapsto\bigcup_{M\in\C}\supp M$$ induces a
bijection between all Serre subcategories of $\mod A$ and all open
subsets of $\Spec^* A$; see \cite{G}. Our second example arises from
Ziegler's work on the model theory of modules \cite{Z}.  The points of
the Ziegler spectrum of $A$ are the isomorphism classes of
indecomposable pure-injective $A$-modules and the closed subsets
correspond to complete theories of modules. The indecomposable
injective modules form a closed subset $\Inj A$, and the assignment
$$\Spec^* A\ni\fp\mapsto E(A/\fp)=\text{ injective envelope of
}A/\fp$$ induces a bijection between all closed subsets of $\Spec^*A$
and all Ziegler closed subsets contained in $\Inj A$; see \cite{P}.
We do not comment any further on the second example, but the first
example is explained in some detail when we discuss the abelian
category of coherent sheaves and the triangulated category of perfect
complexes on a scheme.

Now let us give a brief outline of the contents of this paper.  At the
beginning we introduce the notion of an ideal lattice and study its
prime ideal spectrum. A key observation is the functoriality of the
spectrum and that this functor admits an adjoint. We assign to each
ideal its support and can classify ideals in terms of their support.
Examples of ideal lattices arise from tensor categories which are
abelian or triangulated. We provide a systematic treatment of such
tensor categories and discuss a number of examples from algebraic
geometry and modular representation theory. This is inspired by recent
work of Balmer \cite{B}.

To be more specific, let $\C$ be an abelian or triangulated tensor
category with a tensor identity. We consider the ideal lattice of
thick tensor ideals of $\C$ and its prime ideal spectrum
$\Spec\C$. This space comes naturally equipped with a sheaf of rings
$\OO_\C$ and we can describe the ringed space $(\Spec\C,\OO_\C)$ in
some interesting examples.  An important application says that every
quasi-compact and quasi-separated scheme $X$ can be reconstructed from
the triangulated tensor category of perfect complexes on $X$. This is
a slight generalization of a result of Balmer \cite{B} and based on
the fundamental work of Thomason \cite{TT,T}. On the other hand, it is
the analogue -- with almost identical proof -- of the fact that a
noetherian scheme $X$ can be reconstructed from the abelian tensor
category of coherent sheaves on $X$.

There are a number of interesting examples of triangulated tensor
categories $\C$ where $(\Spec \C,\OO_\C)$ is actually a projective
scheme. We provide a general criterion which explains those
examples. For instance, this result establishes for a finite group $G$
and a field $k$ a conceptual link between the identical
classifications of
\begin{enumerate}
\item[--]thick tensor ideals of the category of perfect complexes over
the graded commutative cohomology ring $H^*(G,k)$, due to Hopkins and
Neeman \cite{H,N}, and
\item[--] thick tensor ideals of the stable category of finite
dimensional $k$-linear representations of $G$, due to Benson, Carlson,
and Rickard \cite{BCR}.
\end{enumerate}
This generalizes to finite group schemes, by recent work of
Friedlander and Pevtsova \cite{FP}.

Our personal motivation for this project stems from the work on
support varieties in non-commutative settings, for instance for
modular representations of finite dimensional algebras. It turns out
that most parts of our theory do not require any commutativity
assumptions. However, the product formula
$$\supp(ab)=\supp(a)\cap\supp(b)=\supp(ba)$$ for the support of two
ideals $a,b$ shows that commutativity is inherent to the
subject, even though we allow $ab\neq ba$.

\subsection*{Acknowldegements} 
The authors are grateful to Paul Balmer for various helpful comments
on this work.

\section{The prime spectrum of an ideal lattice}

\subsection{Ideal lattices}
In this section we introduce the notion of an ideal lattice. The
collection of ideals of some fixed algebraic structure is usually
equipped with two additional structures. We consider the partial
ordering by inclusion and the internal multiplication.

\begin{defn} 
An {\em ideal lattice} is by definition a partially ordered set
$L=(L,\leq)$, together with an associative multiplication $L\times
L\to L$, such that the following holds.
\begin{enumerate}
\item[(L1)] The poset $L$ is a {\em complete lattice}, that is, 
$$\bigvee_{a\in A} a=\sup A\quad\text{and}\quad \bigwedge_{a\in A}
a=\inf A$$ exist in $L$ for every subset $A\subseteq L$.
\item[(L2)] The lattice $L$ is {\em compactly generated}, that is,
every element in $L$ is the supremum of compact elements.  (An element
$a\in L$ is {\em compact}, if for all $A\subseteq L$ with $a\leq \sup
A$ there exists some finite $A'\subseteq A$ with $a\leq\sup A'$.)
\item[(L3)] We have for all $a,b,c\in L$
$$a(b\vee c)=ab\vee ac\quad\text{and}\quad (a\vee b)c=ac\vee bc.$$
\item[(L4)] The element $1=\sup L$ is compact, and $1a=a=a1$ for all $a\in L$.
\item[(L5)] The product of two compact elements is again compact.
\end{enumerate}
A {\em morphism} $\p\colon L\to L'$ of ideal lattices is a map satisfying
\begin{gather*}\label{eq:mor}
\p(\bigvee_{a\in A}a)=\bigvee_{a\in A}\p(a)\quad \text{for}\quad
A\subseteq L, \\
\p(1)=1\quad\text{and}\quad\p(ab)=\p(a)\p(b)\quad\text{for}\quad
a,b\in L.\notag
\end{gather*}
\end{defn}

It is useful to think of a poset $L$ as a category $\L$ where the
objects of $\L$ are the elements of $L$ and
$$\Hom_\L(a,b)\neq\emptyset\iff \card\Hom_\L(a,b)= 1\iff a\leq b$$ for
$a,b\in L$. Note that infimum and supremum in $L$ correspond to
product and coproduct, respectively, in $\L$. Thus a compactly
generated complete lattice is precisely a locally finitely presentable
category $\L$ (in the sense of \cite{GU}) satisfying
$\card\Hom_\L(a,b)\leq 1$ for all $a,b\in\L$. Given an ideal lattice
$L$, the multiplication $L\times L\to L$ corresponds to a tensor product
$\L\times \L\to\L$. A morphism $L\to L'$ of ideal lattices corresponds
to a functor $\L\to\L'$ preserving all colimits and the tensor product.

Next observe that an ideal lattice $L$ is essentially determined by
its subset $L^c$ of compact elements. To make this precise, let $K$ be
a poset and suppose that $\sup A$ exists for every finite subset
$A\subseteq K$.  A non-empty subset $I\subseteq K $ is an {\em ideal}
of $K$ if for all $a,b\in K$
\begin{enumerate}
\item $a\leq b$ and $b\in I$ implies $a\in I$, and
\item $a,b\in I$ implies $a\vee b\in I$.
\end{enumerate}
Given $a\in K$, let $I(a)=\{x\in K\mid x\leq a\}$ denote the {\em
principal ideal} generated by $a$. The set $\widehat K$ of all ideals
of $K$ is called the {\em completion} of $K$. This set is partially
ordered by inclusion and in fact a compactly generated complete
lattice. The map $K\to\widehat K$ sending $a\in K$ to $I(a)$
identifies $K$ with ${\widehat K}^c$.

\begin{lem}
Let $L$ be a compactly generated complete lattice. Then the map
$$L\lto \widehat{L^c},\quad a \mapsto I(a)\cap L^c=\{x\in L\mid x\leq
a\text{ and }x\text{ compact}\},$$ is a lattice isomorphism.
\end{lem}
\begin{proof}
The inverse map sends an ideal $I\in \widehat{L^c}$ to $\sup I$ in $L$.
\end{proof}
We note some immediate consequences which we use frequently without
further reference. Given elements $a,b$ in a compactly generated
complete lattice, we have
\begin{gather}\label{eq:latt}
a=\bigvee_{\genfrac{}{}{0pt}{}{ a'\leq a}{ a' \text{
compact}}}a',\quad \text{and}\\
a\leq b\quad\iff \quad a'\leq b\text{\ \ for all compact\ \ }a'\leq a.\notag
\end{gather}

Now let $L$ be an ideal lattice $L$. The multiplication $L\times L\to L$
restricts to a multiplication $L^c\times L^c\to L^c$. It turns out
that all relevant structure is determined by the multiplication of
compact elements. In our applications, we always have for $a,b\in L$
that $ab=\sup a'b'$ where $a'\leq a$ and $b'\leq b$ run through all
compact elements.

\subsection{The prime spectrum}
Let $L$ be an ideal lattice. We define the spectrum of prime
elements in $L$ and discuss some of its basic properties. An element
$p\neq 1$ in $L$ is called {\em prime} if $ab\leq p$ implies $a\leq p$
or $b\leq p$ for all $a,b\in L$. A subset $S\subseteq L$ is {\em
multiplicative} if $ab\in S$ for all $a,b\in S$. We collect some
elementary facts.

\begin{lem}
An element $p\neq 1$ in $L$ is  prime if and only if $ab\leq p$ implies
$a\leq p$ or $b\leq p$ for all compact $a,b\in L$.
\end{lem}
\begin{proof}
Use (\ref{eq:latt}).
\end{proof}

\begin{lem}\label{le:p-avoid} 
Let $a\in L$ and $S\subseteq L$ be a non-empty multiplicative set of
compact elements. Suppose that $s\not\leq a$ for all $s\in S$. Then there
exists a prime $p\in L$ such that $a\leq p$ and $s\not\leq p$ for all
$s\in S$.
\end{lem}
\begin{proof} Consider the set $A$ of all elements $x\in L$ such that
$a\leq x$ and $s\not\leq x$ for all $s\in S$. The set $A$ is non-empty
and for every chain $B\subseteq A$, we have $\sup B$ in $A$
since the elements in $S$ are compact. Thus $A$ has a maximal element
$p$, by Zorn's lemma. We claim that $p$ is prime. Let $x,x'\in L$ with
$xx'\leq p$. Suppose that $x\not\leq p$ and $x'\not\leq p$. Then we have
$s,s'\in S$ such that $s\leq p\vee x$ and $s'\leq p\vee x'$, by the
maximality of $p$. Therefore
$$ss'\leq (p\vee x)(p\vee x')=pp\vee px'\vee xp\vee xx'\leq p$$ which
contradicts the fact that $p\in A$. Thus $x\leq p$ or $x'\leq p$, and
therefore $p$ is prime.
\end{proof}

An element $a\in L$ is {\em semi-prime} if $bb\leq a$ implies $b\leq
a$ for all $b\in L$.

\begin{lem}\label{le:semiprime}
An element $a\in L$ is semi-prime if and only if $a=\inf V$
for some set $V\subseteq L$ of prime elements.
\end{lem}
\begin{proof} 
Suppose that $a$ is semi-prime and let $V=\{p\in L\mid\mbox{$a\leq p$
and $p$ prime}\}$. For any compact $b\in L$ such that $b\not\leq a$,
consider the multiplicative set $\{b^n\mid n\geq 1\}$. It follows from
Lemma~\ref{le:p-avoid} that there is a prime $p\in V$ such that
$b\not\leq p$.  Thus $a=\inf V$. The other implication is clear.
\end{proof}

We denote by $\Spec L$ the set of prime elements in $L$ and define for
each $a\in L$
$$V(a)=\{p\in\Spec L\mid a\leq p\}\quad\text{and}\quad
D(a)=\{p\in\Spec L\mid a\not\leq p\}.$$ 
The subsets of $\Spec L$ of
the form $V(a)$  are closed under forming
arbitrary intersections and finite unions.  More precisely,
$$V(\bigvee_{i\in\Omega} a_i)=\bigcap_{i\in\Omega}
V(a_i)\quad\text{and}\quad V(ab)=V(a)\cup V(b).$$ Thus we obtain the
{\em Zariski topology} on $\Spec L$ by declaring a subset of $\Spec L$
to be {\em closed} if it is of the form $V(a)$ for some $a\in L$.  The
set $\Spec L$ endowed with this topology is called the {\em prime
spectrum} of $L$.  Note that the sets of the form $D(a)$ with compact
$a\in L$ form a basis of open sets. This is a consequence of the
following lemma.

\begin{lem}\label{le:U(a)}
For $a\in L$, we have
$$V(a)=\bigcap_{\genfrac{}{}{0pt}{}{ b\leq a}{ b \text{
compact}}}V(b)\quad\text{and}\quad D(a)=\bigcup_{\genfrac{}{}{0pt}{}{
b\leq a}{ b \text{ compact}}}D(b).$$
\end{lem}
\begin{proof} 
Use (\ref{eq:latt}).
\end{proof}

\begin{prop}\label{pr:bijection}
The assignments
$$L\ni a\mapsto V(a)=\{p\in\Spec L\mid a\leq p\}
\quad\text{and}\quad \Spec L\supseteq Y\mapsto
\inf Y$$ induce mutually inverse and
order reversing bijections between 
\begin{enumerate} 
\item the set of all semi-prime elements in $L$, and 
\item the set of all closed subsets of $\Spec L$.
\end{enumerate}
\end{prop}
\begin{proof} 
Both maps are well-defined by Lemma~\ref{le:semiprime}.  Given a
semi-prime $a\in L$, the equality $\inf V(a)=a$ is clear since $a$ is
a join of prime elements, by Lemma~\ref{le:semiprime}. Now let
$Y\subseteq \Spec L$. The inclusion $Y\subseteq V(\inf Y)$ is purely
formal. Suppose that $Y$ is of the form $Y=V(a)$ for some $a\in L$.
If $p\in V(\inf Y)$, then $a\leq \inf Y\leq p$ and therefore $p\in
Y$. Thus the proof is complete.
\end{proof}

\begin{cor}\label{co:bijection}
The assignments
$$L\ni a\mapsto D(a)=\bigcup_{\genfrac{}{}{0pt}{}{ b\leq a}{b \text{
compact}}} D(b) \quad\text{and}\quad \Spec L\supseteq Y\mapsto
\bigvee_{\genfrac{}{}{0pt}{}{D(b)\subseteq Y}{b \text{ compact}}}b$$
induce mutually inverse and order preserving bijections between
\begin{enumerate} 
\item the set of all semi-prime elements in $L$, and 
\item the set of all open subsets of $\Spec L$.
\end{enumerate}
\end{cor}
\begin{proof}
We apply Proposition~\ref{pr:bijection} and need to check that for
$V=\Spec L\setminus Y$ and $a\in L$, we have $a\leq \inf V$ if and
only if $D(a)\subseteq Y$. This is clear since $a\leq p$ for all $p\in
V$ is equivalent to $a\not\leq q$ implies $q\in Y$.
\end{proof}

\begin{rem} 
Let $a\in L$ and denote by $\sqrt a=\inf V(a)$ the smallest semi-prime
in $L$ containing $a$. Then we have 
$$\sqrt{ab}=\inf(V(a)\cup V(b))=\sqrt{ba}\quad\text{for}\quad a,b\in
L,$$ even though we do not assume commutativity of the multiplication
in $L$.
\end{rem}

\section{The prime spectrum is spectral}

Recall from \cite{Ho} that a topological space is {\em spectral} if it
is $T_0$ and quasi-compact, the quasi-compact open subsets are closed
under finite intersections and form an open basis, and every non-empty
irreducible closed subset has a generic point. We have the following
basic property of a spectral space.

\begin{lem}\label{le:spectral}
Let $X$ be a spectral space. Endow the underlying set with a new
topology by taking as open sets those of the form
$Y=\bigcup_{i\in\Omega}Y_i$ with quasi-compact open complement
$X\setminus Y_i$ for all $i\in\Omega$, and denote the new space by
$X^*$. Then $X^*$ is spectral and $(X^*)^*=X$.
\end{lem}
\begin{proof} See \cite[Prop.~8]{Ho}.
\end{proof}

Let $L$ be an ideal lattice. Next we show that the space $\Spec L$ is
spectral. We proceed in several steps.

\begin{lem}\label{le:quasicompact}
An open subset of $\Spec L$ is quasi-compact if and only if it is of
the form $D(c)$ for some compact $c\in L$.
\end{lem}
\begin{proof}
Fix an open subset $D(a)$ of $\Spec L$.  Suppose first that $D(a)$ is
quasi-compact.  We have $$D(a)=\bigcup_{\genfrac{}{}{0pt}{}{ b\leq
a}{b \text{ compact}}}D(b)$$ by Lemma~\ref{le:U(a)}, and therefore
$D(a)=D(b)$ for some compact $b\in L$.

Now suppose that $a\in L$ is compact and $D(a)\subseteq \bigcup_{b\in
B}D(b)$ for some subset $B\subseteq L$.  We write $\bar
b=\bigvee_{b\in B}b$ and have $D(a)\subseteq D(\bar b)$.  It follows
from Lemma~\ref{le:p-avoid} that $a^n\leq \bar b$ for some $n\geq
1$. Thus $a^n\leq \bigvee_{b\in B'}b$ for some finite subset $B'\subseteq B$
since $a^n$ is compact. This implies $D(a)\subseteq
\bigcup_{b\in B'}D(b)$, and therefore $D(a)$ is quasi-compact.
\end{proof}

\begin{lem}\label{le:closure}
Let $p,q\in\Spec L$. Then $\overline{\{p\}}=V(p)$. In particular,
$\overline{\{p\}}=\overline{\{q\}}$ implies $p=q$. 
\end{lem}
\begin{proof}
Clear.
\end{proof}

\begin{lem}\label{le:closed}
Let $Y\subseteq\Spec L$ be a non-empty closed subset. 
If $Y$ is irreducible, then $\inf Y$ is prime and $Y=\overline{\{\inf Y\}}$.
\end{lem}
\begin{proof} 
First observe that we have $c\leq \inf Y$ for $c\in L$ if and only if
$Y\subseteq V(c)$. To show that $\inf Y$ is prime, let $ab\leq \inf
Y$. Then $$Y\subseteq V(ab)=V(a)\cup V(b),$$ and we have $Y\subseteq
V(a)$ or $Y\subseteq V(b)$ since $Y$ is irreducible. Thus $a\leq \inf
Y$ or $b\leq \inf Y$. Let $p=\inf Y$ and $Y=V(a)$. Then we have by
construction $a\leq p$ and $Y\subseteq V(p)$. Thus $Y= V(p)=\overline{\{p\}}$,
by Lemma~\ref{le:closure}.
\end{proof}

\begin{prop}\label{pr:spectral}
The prime spectrum $\Spec L$ of an ideal lattice $L$ is spectral.
\end{prop}
\begin{proof} 
The space $\Spec L$ is $T_0$ by Lemma~\ref{le:closure}.  An open
subset of $\Spec L$ is quasi-compact if and only if it is of the form
$D(a)$ for some compact $a\in L$, by Lemma~\ref{le:quasicompact}. Thus
the definition of the topology implies that 
the quasi-compact open subsets  form an open basis which is closed under
finite intersections. Moreover, $\Spec L$ is
quasi-compact since $1=\sup L$ is compact.  If $Y$ is a
non-empty irreducible closed subset, then $Y=\overline{\{p\}}$ for
$p=\inf Y$ by Lemma~\ref{le:closed}.
\end{proof}

There is a close relation between spectral spaces and ideal lattices,
and we make this more precise.  Given a topological space $X$, we
denote by $L_\open(X)$ the lattice of open subsets of $X$ and consider
the multiplication map
$$L_\open(X)\times L_\open(X)\lto L_\open(X),\quad (U,V)\mapsto
UV=U\cap V.$$ Note that the lattice $L_\open(X)$ is complete.

\begin{lem}\label{le:top}
Let $X$ be a space and $U\in L=L_\open(X)$. Then
\begin{enumerate}
\item $U$ is prime in $L$ if and only if $X\setminus U$ is irreducible, and
\item $U$ is compact in $L$ if and only if $U$ is quasi-compact.
\end{enumerate}
\end{lem}
\begin{proof} 
Clear.
\end{proof}

\begin{prop}\label{pr:openlattice}
Let $X$ be a spectral space. Then $L_\open(X)$ is an ideal lattice and
every ideal in $L_\open(X)$ is semi-prime. Moreover, the map
$$X\lto\Spec L_\open(X),\quad x\mapsto X\setminus \overline{\{x\}},$$  
is a homeomorphism.
\end{prop}
\begin{proof}
Using Lemma~\ref{le:top}, the properties of a spectral space can be
translated into the defining properties of an ideal lattice.  Clearly,
every ideal is semi-prime, since $UU=U$ for all $U\in L_\open(X)$. It
is straightforward to check that the given map is a homeomorphism.
\end{proof}

\begin{exm}
Let $A$ be a commutative ring. Then the lattice $L_\id (A)$ of ideals
of $A$ is an ideal lattice and therefore $\Spec A=\Spec L_\id (A)$ is
spectral. More generally, if $X$ is a quasi-compact and
quasi-separated scheme, then the underlying space of $X$ is spectral.
\end{exm}

\section{An adjoint of the functor $\Spec$}\label{se:adjoint}

The prime spectrum of an ideal lattice satisfies a universal property
which we discuss in this section. Then we view the assignment $L\mapsto
\Spec L$ as a functor from ideal lattices to spectral topological
spaces and study its adjoint.

\begin{defn}
A {\em spectrum} of an ideal lattice $L$ is a pair $(X,\d)$ where
$X$ is a topological space and $\d$ is a map which assigns to each
$a\in L$ an open subset $\d(a)\subseteq X$, such that 
\begin{gather}\label{eq:spe}
\d(\bigvee_{a\in A}a)=\bigcup_{a\in A}\d(a)\quad \text{for}\quad
A\subseteq L, \\
\d(1)=X\quad\text{and}\quad\d(ab)=\d(a)\cap\d(b)\quad\text{for}\quad
a,b\in L.\notag
\end{gather}
A {\em morphism} $f\colon (X,\d)\to (X',\d')$ of
spectra is a continuous map $f\colon X\to X'$ such that
$\d(a)=f^{-1}(\d'(a))$ for all $a\in L$. Such a morphism is an
isomorphism if and only if $f\colon X\to X'$ is a homeomorphism.
\end{defn}

\begin{thm}\label{th:specfinal}
Let $L$ be an ideal lattice.  Then the pair $(\Spec L,D)$ is a
spectrum of $L$.  For every spectrum $(X,\d)$ of $L$, there exists a
unique continuous map $f\colon X\to\Spec L$ such that
$\d(a)=f^{-1}(D(a))$ for every $a\in L$. The map $f$ is defined by
$$f(x)=\bigvee_{\genfrac{}{}{0pt}{}{x\not\in\d(c)}{c \text{
compact}}}c \quad\textrm{for}\quad x\in X.$$
\end{thm}
\begin{proof}
Clearly, the pair $(\Spec L,D)$ is a spectrum. Now let $(X,\d)$ be a
spectrum of $L$. We show that for each $x\in X$ the element
$f(x)=\bigvee_{x\not\in\d(c)}c$ is prime.  First observe that
$f(x)\neq 1$ since $\d(1)=X$.  Suppose that $ab\leq f(x)$, and
we may assume that $a,b$ are compact.  Using (\ref{eq:spe}), observe
that the compact $c\in L$ with $x\not\in\d(c)$ form a directed set.
Thus $ab\leq c$ for some compact $c\in L$ with $x\not\in\d(c)$. Using
again (\ref{eq:spe}), we have $$\d(a)\cap\d(b)= \d(ab)\subseteq\d(c)$$
and therefore $x\not\in \d(a)$ or $x\not\in \d(b)$. We conclude that
$a\leq f(x)$ or $b\leq f(x)$. The definition of $f$ implies
$\d(a)=f^{-1}(D(a))$ for every $a\in L$, since
$$ x\in\d(a)\iff a\not\leq f(x)\iff f(x)\in D(a) \iff x\in
f^{-1}(D(a)).$$ The continuity of $f$ follows from the fact that the
sets $D(a)$ with $a\in L$ are precisely the open subsets of $\Spec
L$. Now let $f_1,f_2\colon X\to\Spec L$ be two maps satisfying
$$f_1^{-1}(D(a))=\d(a)=f_2^{-1}(D(a))$$ for every $a\in L$. Fix $x\in
X$. Then we have
$$\overline{\{f_1(x)\}}=\bigcap_{f_1(x)\not\in D(a)}V(a)=
\bigcap_{f_2(x)\not\in D(a)}V(a)=\overline{\{f_2(x)\}}.$$ This implies
$f_1(x)=f_2(x)$ by Lemma~\ref{le:closure}, and therefore the proof is
complete.
\end{proof}

Our first application says that $\Spec$ is actually a functor into the
category of topological spaces.

\begin{lem}
A morphism of ideal lattices $\p\colon L\to L'$ induces a unique
continuous map $\Spec\p\colon\Spec L'\to\Spec L$ such that
$$D(\p(a))=(\Spec\p)^{-1}D(a)\quad\text{for}\quad a\in L.$$
\end{lem}
\begin{proof}
The pair $(\Spec L',D\comp\p)$ is a spectrum of $L$. Now apply 
Theorem~\ref{th:specfinal}. Note that we can compute more explicitly
$$(\Spec\p)p=\sup\{a\in L\mid \p(a)\leq p\}\quad\text{for}\quad p\in
\Spec L'.$$
\end{proof}

The universal property of the prime spectrum yields an adjoint functor
for $\Spec$.

\begin{thm}\label{th:adj}
We have an adjoint pair of contravariant functors
\[
\xymatrix{\mathbf{Lat}_\mathbf{id}
  \ar@<1ex>[rr]^-{\Spec}&&\mathbf{Top}_\mathbf{sp}
  \ar@<1ex>[ll]^-{L_\open}}
\]
between the category of ideal lattices and the category of spectral
spaces. More precisely, for an ideal lattice $L$ and a spectral space
$X$, there are mutually inverse bijections
\[
\xymatrix{\Hom_{\mathbf{Lat}_\mathbf{id}}(L,L_\open(X))
  \ar@<1ex>[rr]^-{\Si}&&\Hom_{\mathbf{Top}_\mathbf{sp}}(X,\Spec L).
  \ar@<1ex>[ll]^-{\La}}
\]
The functor $L_\open$ is fully faithful, and an ideal lattice $L$ is
isomorphic to one of the form $L_\open(X)$ if and only if every ideal
in $L$ is semi-prime.
\end{thm}
\begin{proof}
Both functors are well-defined by Propositions~\ref{pr:spectral} and
\ref{pr:openlattice}. The maps $\Si$ and $\La$ are defined by 
$$(\Si\p)(x)=\bigvee_{\genfrac{}{}{0pt}{}{x\not\in \p(c)}{c \text{
compact}}}c\quad\text{and}\quad (\La f)(a)= f^{-1}(D(a))
\quad\text{for}\quad x\in X, a\in L.$$ It follows from
Theorem~\ref{th:specfinal} that both maps are mutually inverse
bijections.  Next observe that $L_\open$ is fully faithful. This
follows from the fact that the adjunction morphism $X\to\Spec
L_\open(X)$ is a homeomorphism; see also
Proposition~\ref{pr:openlattice}. It remains to describe the image of
$L_\open$. Clearly, every ideal in $L_\open (X)$ is semi-prime; see
Proposition~\ref{pr:openlattice}. Conversely, if every ideal in $L$ is
semi-prime, then the adjunction morphism $L\to L_\open(\Spec L)$ is an
isomorphism, by Corollary~\ref{co:bijection}. Thus the proof is
complete.
\end{proof}

We observe that a morphism $\p\colon L\to L_\open(X)$ is determined by
its restriction to the subset $L^c$ of compact elements since
$$ \p(a)=\bigcup_{\genfrac{}{}{0pt}{}{ b\leq a}{ b \text{
compact}}}\p(b)\quad\text{for}\quad a\in L.$$ In particular, we obtain
an equivalent definition of a spectrum of $L$ by taking maps $\d\colon
L^c\to L_\open (X)$ satisfying
\begin{gather*}
\d(a\vee b)=\d(a)\cup\d(b)\quad \text{for}\quad
a,b\in L, \\
\d(1)=X\quad\text{and}\quad\d(ab)=\d(a)\cap\d(b)\quad\text{for}\quad
a,b\in L.\notag
\end{gather*}
This observation leads naturally to the concept of a support datum.

\section{Support data}

Let $L$ be an ideal lattice. We have seen that the space $\Spec L$ is
spectral and, in view of our applications, we consider from now on the
`opposite' topology on $X=\Spec L$. To be precise, we let
$\Spec^* L=X^*$ where the points of $X^*$ and $X$ coincide and
$Y\subseteq X^*$ is by definition {\em open} if
$Y=\bigcup_{i\in\Om}Y_i$ with quasi-compact Zariski-open complement
$X\setminus Y_i$ for all $i\in\Om$; see Lemma~\ref{le:spectral}. For
$a\in L$, we call
$$\supp (a)=\{p\in\Spec^* L\mid a\not\leq p\}$$ the {\em support} of
$a$ and observe that $\supp (a)$ is closed if $a$ is compact. Let us
reformulate the classification of semi-prime ideals in terms of the
topology $\Spec^* L$.

\begin{prop}\label{pr:bijection*}
The assignments
$$L\ni a\mapsto \supp(a)=\bigcup_{\genfrac{}{}{0pt}{}{ b\leq a}{ b
\text{ compact}}} \supp(b) \quad\text{and}\quad \Spec^* L\supseteq
Y\mapsto \bigvee_{\genfrac{}{}{0pt}{}{\supp(b)\subseteq Y}{b \text{
compact}}}b$$ induce mutually inverse and order preserving
bijections between
\begin{enumerate} 
\item the set of all semi-prime elements in $L$, and 
\item the set of all subsets $Y\subseteq \Spec^* L$ of the form
$Y=\bigcup_{i\in \Omega} Y_i$ with quasi-compact open complement
$\Spec^* L\setminus Y_i$ for all $i\in \Omega$.
\end{enumerate}
\end{prop}
\begin{proof}
Use Corollary~\ref{co:bijection} and observe that the subsets
$Y=\bigcup_{i\in \Omega} Y_i$ with quasi-compact open complement
$\Spec^* L\setminus Y_i$ are precisely the open subsets of $\Spec
L=(\Spec^*  L)^*$, by Lemma~\ref{le:spectral}.
\end{proof}

Next we introduce for an ideal lattice the concept of a support
datum. This is inspired by Balmer's definition of a support datum on a
triangulated tensor category \cite[Defn.~3.1]{B}. The subsequent
theorem is the analogue of \cite[Thm.~3.2]{B}.

\begin{defn}
A {\em support datum} on an ideal lattice $L$ is a pair $(X,\s)$ where
$X$ is a topological space and $\s$ is a map which assigns to each
compact $a\in L$ a closed subset $\s(a)\subseteq X$, such that 
\begin{gather}\label{eq:supp}
\s(a\vee b)=\s(a)\cup\s(b)\quad \text{for}\quad a,b\in L, \\
\s(1)=X\quad\text{and}\quad\s(ab)=\s(a)\cap\s(b)\quad \text{for}\quad
a,b\in L.\notag
\end{gather}
A {\em morphism} $f\colon (X,\s)\to (X',\s')$ of
support data is a continuous map $f\colon X\to X'$ such that
$\s(a)=f^{-1}(\s'(a))$ for all compact $a\in L$. Such a morphism is an
isomorphism if and only if $f\colon X\to X'$ is a homeomorphism.
\end{defn}

The following result complements the universal property of the pair $(\Spec
L,D)$ of Theorem~\ref{th:specfinal} and the proof is almost the same.

\begin{thm}\label{th:final}
Let $L$ be an ideal lattice. 
Then the pair $(\Spec^*  L,\supp)$ is a support datum on $L$.  For every
support datum $(X,\s)$ on $L$, there exists a unique continuous map
$f\colon X\to\Spec^*  L$ such that $\s(a)=f^{-1}(\supp(a))$ for every
compact $a\in L$. The map $f$ is defined by
$$f(x)=\bigvee_{\genfrac{}{}{0pt}{}{x\not\in\s(c)}{c \text{
compact}}}c \quad\textrm{for}\quad x\in X.$$
\end{thm}
\begin{proof}
Clearly, the pair $(\Spec^* L,\supp)$ is a support datum. Now let
$(X,\s)$ be a support datum on $L$. We show that for each $x\in X$ the
element $f(x)=\bigvee_{x\not\in\s(c)}c$ is prime.  First observe that
$f(x)\neq 1$ since $\s(1)=X$. Suppose that $ab\leq f(x)$, and we may
assume that $a,b$ are compact.  Using (\ref{eq:supp}), observe that
the compact $c\in L$ with $x\not\in\s(c)$ form a directed set.
Thus $ab\leq c$ for some compact $c\in L$ with $x\not\in\s(c)$. Using
again (\ref{eq:supp}), we have $$\s(a)\cap\s(b)=
\s(ab)\subseteq\s(c)$$ and therefore $x\not\in \s(a)$ or $x\not\in
\s(b)$. We conclude that $a\leq f(x)$ or $b\leq f(x)$. The definition
of $f$ implies $\s(a)=f^{-1}(\supp(a))$ for every compact $a\in L$,
since
$$ x\in\s(a)\iff a\not\leq f(x)\iff f(x)\in \supp(a) \iff x\in
f^{-1}(\supp(a)).$$ The continuity of $f$ follows from the fact that the
sets $\supp (a)$ with compact $a\in L$ form a basis of closed sets for
the topology on $\Spec^* L$. Now let $f_1,f_2\colon X\to\Spec^* L$ be
two maps satisfying
$$f_1^{-1}(\supp(a))=\s(a)=f_2^{-1}(\supp(a))$$ for every compact
$a\in L$. Fix $x\in X$. Then we have
$$\overline{\{f_1(x)\}}=\bigcap_{\genfrac{}{}{0pt}{}{f_1(x)\in\supp(a)}{
a \text{ compact}}}\supp(a)=
\bigcap_{\genfrac{}{}{0pt}{}{f_2(x)\in\supp(a)}{ a \text{
compact}}}\supp(a)=\overline{\{f_2(x)\}}.$$ This implies
$f_1(x)=f_2(x)$ since the space $\Spec^*L$ is $T_0$, by
Proposition~\ref{pr:spectral} and Lemma~\ref{le:spectral}
\end{proof}

\section{Classifying support data}

Let $L$ be an ideal lattice. A support datum $(X,\s)$ on $L$ is called
{\em classifying} if the space $X$ is spectral and the assignments
$$L\ni a\mapsto \bigcup_{\genfrac{}{}{0pt}{}{ b\leq a}{ b \text{
compact}}}\s(b)\quad\text{and}\quad X\supseteq Y\mapsto
\bigvee_{\genfrac{}{}{0pt}{}{\s(b)\subseteq Y}{b \text{ compact}}}b$$
induce bijections between
\begin{enumerate} 
\item the set of all semi-prime elements in $L$, and
\item the set of all subsets $Y\subseteq X$ of the form
$Y=\bigcup_{i\in \Omega} Y_i$ with quasi-compact open complement
$X\setminus Y_i$ for all $i\in \Omega$.
\end{enumerate}

\begin{prop}
Let $f\colon (X,\s)\to (X',\s')$ be a morphism of support
data. If both support data are classifying, then the map $f\colon X\to
X'$ is a homeomorphism.
\end{prop}
\begin{proof}
Let $Y\subseteq X$ and $Y'\subseteq X'$ be subsets which are unions of
subsets with quasi-compact open complement, and suppose
$$\bigvee_{\genfrac{}{}{0pt}{}{\s(b)\subseteq Y}{b \text{
compact}}}b=a= \bigvee_{\genfrac{}{}{0pt}{}{\s'(b)\subseteq Y'}{b
\text{ compact}}}b.$$ Then we have
$$Y= \bigcup_{\genfrac{}{}{0pt}{}{ b\leq a}{ b \text{
compact}}}\s(b)\quad\text{and}\quad Y'= \bigcup_{\genfrac{}{}{0pt}{}{
b\leq a}{ b \text{ compact}}}\s'(b).$$ This implies
$$f^{-1}(Y')=\bigcup_{\genfrac{}{}{0pt}{}{ b\leq a}{ b \text{
compact}}} f^{-1}(\s'(b))=\bigcup_{\genfrac{}{}{0pt}{}{ b\leq a}{ b
\text{ compact}}}\s(b)=Y.$$ It follows that the map $Y\mapsto
f^{-1}(Y)$ induces an inclusion preserving bijection between the open
subsets of $X^*$ and $(X')^*$.  In fact, we use that $X,X'$ are
spectral and apply Lemma~\ref{le:spectral}. Thus $f$ is a
homeomorphism $X^*\to (X')^*$ and therefore also a homeomorphism $X\to
X'$.
\end{proof}

The following consequence is the analogue of \cite[Thm.~5.2]{B}. Note
that we do not assume that the support space is noetherian.

\begin{cor}\label{co:class}
A support datum $(X,\s)$ on $L$ is classifying if and only if the
canonical morphism $(X,\s)\to (\Spec^*  L,\supp)$ is an isomorphism.
\end{cor}

\section{Thick tensor ideals}

In this section we consider an additive category with a tensor product
and study its collection of ideals. If there is an additional abelian
or triangulated structure, then we consider those tensor ideals which
are in addition thick subcategories. 

In this paper, all categories are assumed to be {\em small}, that is,
the isomorphism classes of objects form a set (in some fixed
universe).

\subsection{Sublattices of an ideal lattice}

Let $L$ be an ideal lattice. We fix a subset $L'\subseteq L$ satisfying the
following conditions.
\begin{enumerate}
\item[(L$\wedge$)] If $A\subseteq L'$, then $\inf A\in L'$.
\item[(L$\vee$)] If $A\subseteq L'$ is directed, then $\sup A\in L'$.
\end{enumerate}
We consider on $L'$ the partial order induced from the partial order
on $L$ and define the map
$$\pi\colon L\lto L',\quad a\mapsto\bigwedge_{a\leq a'\in L'}a'.$$
Note that we have
\begin{equation}\label{eq:L}
\pi(a)\leq a' \iff a\leq a'\quad\text{\ for\ }\quad a\in L,\,a'\in L'.
\end{equation}
Thus $\pi$ is a left adjoint of the inclusion $L'\to L$ if we think of
posets as categories.  Moreover, we have $1=\inf\emptyset\in L'$.

\begin{lem}\label{le:sub1}
The poset $L'$ is a complete and compactly generated lattice. Every
compact element in $L'$ is of the form $\pi(a)$ for some compact $a\in
L$.
\end{lem}
\begin{proof}
Let $A\subseteq L'$. Then we use (L$\wedge$) to compute the infimum
$\inf A$ in $L'$ and have $$\sup A=\inf\{a'\in L'\mid\mbox{$ a\leq a'$
for all $a\in A$}\}.$$ It follows from (\ref{eq:L}) and (L$\vee$) that
$\pi$ preserves compactness and that each element in $L'$ is the
supremum of compact elements. Thus $L'$ is compactly generated.  If
$a'\in L'$ is compact, write $a'=\bigvee_i a_i$ as directed union of
all compact elements $a_i\leq a'$ in $L$ and use that
$a'=\pi(a')=\bigvee_i\pi(a_i)$ equals $\pi(a_i)$ for some index $i$.
\end{proof}

Given $a,b\in L'$, we define their product in $L'$ as
$$a\cdot b=\pi(ab)$$
and use a dot to distinguish it from the product in $L$. We make a
further assumption.
\begin{enumerate}
\item[(L$\pi$)] Given $a,b\in L$, we have
  $\pi(a\pi(b))=\pi(ab)=\pi(\pi(a)b)$. 
\end{enumerate}

\begin{lem}\label{le:sub2}
Let $a,b,c\in L'$. Then we have
\begin{enumerate}
\item $(a\cdot b)\cdot c=a\cdot (b\cdot c)$,
\item $a\cdot 1=a= 1\cdot a$,
\item $a\cdot(b\vee c)=(a\cdot b)\vee (a\cdot c)$ and $(a \vee b)\cdot
c=(a\cdot c)\vee (b\cdot c)$.
\end{enumerate}
\end{lem}
\begin{proof}
Clear.
\end{proof}

\begin{prop}\label{pr:sub}
Let $L$ be an ideal lattice and $L'$ be a subset satisfying the
conditions \emph{(L$\wedge$)}, \emph{(L$\vee$)}, and
\emph{(L$\pi$)}. Then $L'$ inherits from $L$ the structure of an ideal
lattice.
\end{prop}
\begin{proof}
We apply Lemmas~\ref{le:sub1} and \ref{le:sub2}.  Thus $L'$ satisfies
(L1) -- (L4).  To check (L5), let $a',b'\in L'$ be compact and choose
compact elements $a,b\in L$ with $\pi(a)=a'$ and $\pi(b)=b'$.  Now we
obtain
$$a'\cdot b'=\pi(a'b')=\pi(\pi(a)\pi(b))=\pi(ab).$$ The element $ab$
is compact in $L$, and $\pi$ preserves compactness by (L$\vee$). Thus
$a'\cdot b'$ is compact in $L'$.
\end{proof}

\subsection{The ideal lattice of a semi-ring}

Let $A=(A,+,\cdot)$ be a {\em semi-ring}, that is, $A$ is a set
together with two associative binary operations with identities
(denoted by $0$ and $1$) such that the addition is commutative and
distributivity holds.  A subset $I\subseteq A$ containing $0$ is by
definition an {\em ideal} if for all $x,y\in A$
\begin{enumerate}
\item $x\in  I$ and $y\in I$ implies $x+y\in I$, and
\item $x\in  I$ or $y\in I$ implies $xy\in I$.
\end{enumerate}
The ideals of $A$ are partially ordered by inclusion and form a
lattice which we denote by $L_\id(A)$.  Given $ I, J\in L_\id(A)$, we
define
$$ I J=\{\sum_i x_iy_i\mid
x_i\in I,y_i\in J\}\quad\text{and}\quad  I+ J=\{x+y\mid
x\in I,y\in J\}.$$ Note that $ I+ J= I\vee J$. 

\begin{prop}\label{pr:semiring}
Let $A$ be semi-ring. Then the lattice $L_\id(A)$ of ideals satisfies
the conditions \emph{(L1)} -- \emph{(L4)}.  An ideal in $L_\id(A)$ is
compact if and only if it is finitely generated.  If $A$ is
commutative, then condition \emph{(L5)} is satisfied.
\end{prop}
\begin{proof} 
The proof is straightforward. To identify the compact elements, one
uses that $\bigvee_i I_i=\bigcup_i I_i$ for any directed set of
ideals $ I_i$. To show (L5), let $ I=\langle I_0\rangle$ and
$ J=\langle J_0\rangle$ be ideals generated by subsets $ I_0$ and
$ J_0$, respectively. If $A$ is commutative, then $ I J=\langle
xy\mid x\in I_0,y\in J_0\rangle$. Therefore (L5) holds.
\end{proof}

\begin{exm}
Let $A$ be a (not necessarily commutative) ring and suppose $A$
satisfies the ascending chain condition on ideals. Then the
lattice $L_\id(A)$ of ideals is an ideal lattice. Note that we
have the following weak commutativity: $\sqrt{ I J}=\sqrt{ J I}$
for any pair $ I, J$ of ideals.
\end{exm}

\subsection{Thick tensor ideals}

Let $\C=(\C,\otimes,e)$ be an additive category with a {\em tensor
product}.  To be precise, we have an additive bifunctor
$\C\times\C\to\C$ and a natural isomorphism $(x\otimes y)\otimes
z\xto{\sim}x\otimes(y\otimes z)$. In addition, we require the
existence of a {\em tensor identity} $e$, that is, we have natural
isomorphisms $x\otimes e\xto{\sim}x$ and $e\otimes y\xto{\sim}y$
satisfying the Pentagon Axiom and the Triangle Axiom.

Denote by $C$ the set of isomorphism classes of objects in $\C$, and
let $x+y=x\amalg y$ and $xy=x\otimes y$ for $x,y\in C$. Then $C$ is a
semi-ring, and we shall identify $\C$ and $C$ whenever it is
convenient.  A {\em tensor ideal} of $\C$ is a full additive
subcategory $\D$ such that for all $x,y\in\C$, we have $x\otimes
y\in\D$ if $x\in\D$ or $y\in\D$. Note that the tensor ideals in $\C$
are precisely the ideals of the semi-ring $C$. We denote by
$L_\id(\C)$ the lattice of tensor ideals of $\C$ and define the
multiplication of tensor ideals as in $L_\id(C)$.

Now suppose that there exists some additional exact or triangulated
structure on $\C$. A full subcategory of $\C$ is called {\em thick} if
it is `compatible' with this additional structure; see below.  We view
the thick tensor ideals as a subset of $L_\id(\C)$ and denote it by
$L_\thick(\C)$. An ideal $\D\in L_\thick(\C)$ is {\em generated} by a
class $\D_0$ of objects and we write
$$\D=\langle\D_0\rangle$$ if $\D$ is the smallest thick tensor ideal
containing $\D_0$. The product of $\D_1,\D_2$ in $L_\thick(\C)$ is by
definition $\langle\D_1\D_2\rangle$ where $\D_1\D_2$ is computed in
$L_\id(\C)$.

\subsection{Abelian and triangulated tensor categories}

Let $\C$ be an abelian category, or more generally, an exact category
in the sense of Quillen. A full subcategory $\D$ is called {\em
thick} if for every exact sequence $0\to x'\to x\to x''\to 0$ in $\C$,
we have $x\in\D$ if and only if $x',x''\in\D$.  Now suppose that there
is a tensor product $\otimes$ defined on $\C$.

\begin{prop}\label{pr:abel}
Let $\C$ be an abelian category, or more generally an exact
category, with a tensor product which is exact in each variable.
Suppose that either the tensor product is commutative, or that there
exists an object $c\in\C$ such that there is no proper thick
subcategory of $\C$ containing $c$. Then the thick tensor ideals in
$\C$ form an ideal lattice. Moreover, an ideal is compact if and only
if it is generated by a single object.
\end{prop}
\begin{proof}
It is straightforward to check the conditions (L$\wedge$) and
(L$\vee$) for the subset of thick tensor ideals $L_\thick(\C)$ in
$L_\id(\C)$. To verify (L$\pi$), observe that for a tensor ideal
$\D\in L_\id(\C)$, the thick subcategory generated by $\D$ equals
$\pi(\D)$. Here we use the exactness of the tensor product. We deduce
from Propositions~\ref{pr:sub} and \ref{pr:semiring} that
$L_\thick(\C)$ is an ideal lattice. Note that a finitely generated
ideal $\langle x_1,\ldots,x_n\rangle$ is generated by
$x_1\amalg\ldots\amalg x_n$. Thus compact and cyclic ideals coincide.
Finally, if $\otimes$ is not commutative, then (L5) follows from the
identity $\langle x\rangle\langle y\rangle=\langle x\otimes c\otimes
y\rangle$, assuming that $c$ generates $\C$.
\end{proof}

Now let $\C$ be a triangulated category. A full subcategory $\D$ is
called {\em thick} if $\D$ is a triangulated subcategory and for each
$x\in\D$ a decomposition $x=x_1\amalg x_2$ implies $x_1,x_2\in\D$. Now
suppose that there is a tensor product $\otimes$ defined on $\C$.

\begin{prop}\label{pr:tria}
Let $\C$ be a  triangulated category with a tensor product which is
exact in each variable.  Suppose that either the tensor product is
commutative, or that there exists an object $c\in\C$ such that there
is no proper thick subcategory of $\C$ containing $c$. Then the thick
tensor ideals in $\C$ form an ideal lattice. Moreover, an ideal is
compact if and only if it is generated by a single object.
\end{prop}
\begin{proof}
The proof is the same as that of Proposition~\ref{pr:abel}.
\end{proof}

\begin{rem}\label{re:sub}
Let $(\C,\otimes,e)$ be an abelian or triangulated tensor category
with exact tensor product and tensor identity $e$. Let $\E$ be the
thick subcategory generated by $e$. Then we have $x\otimes y\in\E$ for
all $x,y\in\E$ and therefore $(\E,\otimes,e)$ is a category satisfying
the assumptions from Propositions~\ref{pr:abel} or \ref{pr:tria}.
\end{rem}

\subsection{Support data}
Let $\C$ be an abelian or triangulated tensor category.  We assume
from now on that the lattice of thick tensor ideals of $\C$ is an
ideal lattice, for instance by imposing the assumptions from
Propositions~\ref{pr:abel} or \ref{pr:tria}.  We write
$$\Spec \C=\Spec^* L_\thick(\C)$$ for the spectrum of prime ideals.
Note that we keep the `opposite' of the Zariski topology in view of
our applications. This practice is in accordance with Balmer's notion
of a spectrum in \cite{B}.  The compact ideals in $L_\thick(\C)$ are
precisely the ideals $\langle x\rangle$ generated by a single object
$x\in\C$. We write
$$\supp(x)=\supp(\langle x\rangle)=\{\P\in\Spec \C\mid
x\not\in\P\}\quad\text{for}\quad x\in\C$$ and call this subset of
$\Spec\C$ the {\em support} of $x$. It is convenient to work with
support data defined on objects of $\C$ instead of support data
defined on ideals of $\C$. This motivates the following definition
from \cite{B}.

\begin{defn}
A {\em support datum} on $\C$ is a pair $(X,\t)$ where $X$ is a
topological space and $\t$ is a map which assigns to each object
$x\in\C$ a closed subset $\t(x)\subseteq X$, such that for all
$x,y\in\C$
\begin{gather*}
\t(x)=\bigcup_{x'\in\langle x\rangle}\t(x'),\quad \t(x\amalg
y)=\t(x)\cup \t(y),\\ \t(e)=X\quad\text{and}\quad\t(x\otimes y
)=\t(x)\cap\t(y).
\end{gather*}
A {\em morphism} $f\colon (X,\t)\to (X',\t')$ of support data is a
continuous map $f\colon X\to X'$ such that $\t(x)=f^{-1}(\t'(x))$ for
all $x\in\C$.
\end{defn}

\begin{lem}
Let $\C$ be an abelian or triangulated tensor category satisfying the
assumptions from Propositions~\ref{pr:abel} or \ref{pr:tria}.  
\begin{enumerate}
\item If $(X,\t)$ is a support datum on $\C$, then $\s(\langle
x\rangle)=\t(x)$ defines a support datum on the lattice of thick tensor
ideals of $\C$.
\item If $(X,\s)$ is a support datum on the
lattice of thick tensor ideals of $\C$, then $\t(x)=\s(\langle
x\rangle)$ defines a support datum on $\C$.
\end{enumerate}
\end{lem}
\begin{proof}
We start with a support datum $(X,\t)$ on $\C$. The map $\s$ on
compact ideals of $\C$ is well-defined because of the condition
$\t(x)=\bigcup_{x'\in\langle x\rangle}\t(x')$. Now compute for compact
ideals $\langle x\rangle$ and $\langle y\rangle$
$$\s(\langle x\rangle)\cup \s(\langle y\rangle)= \t(x)\cup
\t(y)=\t(x\amalg y)=\s(\langle x\amalg y\rangle)=\s(\langle
x\rangle\vee\langle y\rangle),$$
and if the tensor product is commutative
$$\s(\langle x\rangle)\cap \s(\langle y\rangle)= \t(x)\cap
\t(y)=\t(x\otimes y)=\s(\langle x\otimes y\rangle)=\s(\langle
x\rangle\langle y\rangle),$$ using that $\langle x\rangle\langle
y\rangle=\langle x\otimes y\rangle$.  In the non-commutative case, we
have $\langle x\rangle\langle y\rangle=\langle x\otimes c\otimes
y\rangle$ for some $c\in\C$ with $\langle c\rangle =\C$. Thus
\begin{multline*}
\s(\langle x\rangle)\cap \s(\langle y\rangle)= \t(x)\cap \t(y)=
\t(x)\cap\t(c)\cap \t(y)\\ =\t(x\otimes c\otimes y)=\s(\langle x\otimes
c\otimes y\rangle)=\s(\langle x\rangle\langle y\rangle).
\end{multline*} 
Finally, we have
$$\s(1)=\s(\langle e\rangle)=\t(e)=X.$$ We conclude that $(X,\s)$ is
support datum on the lattice of thick tensor ideals of $\C$. The proof
of the converse is analogous.
\end{proof}

From now on we do not distinguish between support data on the lattice
of thick tensor ideals $L_\thick(\C)$ and support data on $\C$.  We
leave it to the interested reader to reformulate our general results
about ideal lattices for the lattice $L_\thick(\C)$ and its spectrum
$\Spec\C$.

\section{The structure sheaf of a tensor category}

Let $\C$ be an abelian or triangulated tensor category
with tensor identity $e$.  Following \cite{B}, we define a structure
sheaf on $\Spec\C$ as follows.  For an open subset $U\subseteq \Spec
\C$, let
$$\C_U=\{x\in\C\mid\supp(x)\cap U=\emptyset\}$$ and observe that
$\C_U$ is a thick tensor ideal.  We obtain a presheaf of rings on
$\Spec \C$ by
$$U\mapsto \End_{\C/\C_U}(e).$$ If $V\subseteq U$ are open subsets,
then the restriction map $\End_{\C/\C_U}(e)\to \End_{\C/\C_V}(e)$ is
induced by the quotient functor $\C/\C_U\to\C/\C_V$.  The
sheafification is called the {\em structure sheaf} of $\C$ and is
denoted by $\OO_\C$. Note that the endomorphism ring of a tensor
identity is commutative, if the tensor product is commutative, or if
$\C$ is a suspended tensor category; see for instance
\cite[Thm.~1.7]{S}. Next observe that
$$\OO_{\C,\P}\cong\End_{\C/\P}(e) \quad\text{for each}\quad
\P\in\Spec \C.$$ This is an immediate consequence of the following
lemma.

\begin{lem}\label{le:stalk}
Let $\C$ be an abelian or triangulated tensor category and
$\P\in\Spec \C$. Then
$$\li_{\P\in
U}\Hom_{\C/{\C_U}}(x,y)\xto{\sim}\Hom_{\C/{\P}}(x,y)\quad\text{for
all}\quad x,y\in\C,$$ where $U$ runs through all (quasi-compact) open
subsets containing $\P$.
\end{lem}
\begin{proof}
Use that $\P=\bigcup_{\P\in U}\C_U$.
\end{proof}

Now we discuss briefly the functoriality of the spectrum.

\begin{lem}
Let $F\colon \C\to\C'$ be an exact tensor functor. 
\begin{enumerate}
\item $F$ induces a unique continuous map $f\colon\Spec
\C'\to\Spec \C$ such that
$$\supp(Fx)=f^{-1}(\supp(x))\quad\text{for}\quad x\in\C.$$ The map
sends $\P\in\Spec\C'$ to $F^{-1}(\P)$.
\item $F$ induces a morphism of ringed spaces
$$(f,f^\sharp)\colon(\Spec \C',\OO_{\C'})\to(\Spec \C,\OO_\C).$$
\end{enumerate}
\end{lem}
\begin{proof} 
(1) The map sending $x\in\C$ to $\supp (Fx)$ is a support datum on
$\C$. Now apply Theorem~\ref{th:final} to obtain a continuous map
$\Spec \C'\to\Spec \C$.  

(2) Let $U\subseteq \Spec \C$ be open. Then $F$ maps $\C_U$ to
$\C_{f^{-1}U}$ and induces a functor
$\C/\C_U\to\C'/\C'_{f^{-1}U}$. This functor induces a homomorphism
$\End_{\C/\C_U}(e)\to\End_{\C'/\C'_{f^{-1}U}}(e')$. Thus we obtain a
 morphism $f^\sharp\colon \OO_\C\to f_*\OO_{\C'}$. 
\end{proof}

\section{Applications to schemes}

\subsection{Coherent sheaves on a scheme}

We consider a noetherian scheme $X$ and reconstruct it from the
abelian tensor category $\coh X$ of coherent $\OO_X$-modules. This is
based on the following well-known classification of all thick
subcategories of $\coh X$. Given $x\in\coh X$, we write
$$\supp_X(x)=\{P\in X\mid x_P\neq 0\}.$$

\begin{prop}\label{pr:cohX}
Let $X$ be a noetherian scheme. The assignments
$$\coh X\supseteq
\D\mapsto\bigcup_{x\in\D}\supp_X(x)\quad\text{and}\quad X\supseteq
Y\mapsto \{x\in\coh X\mid\supp_X(x)\subseteq Y\}$$ induce bijections
between
\begin{enumerate} 
\item the set of all thick subcategories of $\coh X$, and
\item the set of all subsets $Y\subseteq X$ of the form
$Y=\bigcup_{i\in \Omega} Y_i$ with quasi-compact open complement
$X\setminus Y_i$ for all $i\in \Omega$.
\end{enumerate}
\end{prop}
\begin{proof}
See \cite[Prop.~VI.4]{G}.
\end{proof}

Note that every open subset of a noetherian space is quasi-compact.
Nonetheless, the above formulation is appropriate because it
generalizes to schemes which are not necessarily noetherian.

The abelian category $\coh X$ carries a commutative tensor product
$\otimes_{\OO_X}$, and we deduce from the classification of thick
subcategories the following properties.

\begin{prop}
Let $X$ be a noetherian scheme and $\C=\coh X$. Then every thick
subcategory of $\C$ is a tensor ideal and the thick tensor ideals of $\C$
form an ideal lattice.
\end{prop}
\begin{proof}
We apply Proposition~\ref{pr:cohX}.  The formula
$$\supp_X(x\otimes_{\OO_X}y)=\supp_X(x)\cap\supp_X(y)$$ shows that
every thick subcategory is a tensor ideal. The space $X$ is spectral
because the scheme is noetherian. Thus $X^*$ is spectral and
$L_\open(X^*)$ is an ideal lattice, by
Proposition~\ref{pr:openlattice}.  We have an isomorphism $L_\open
(X^*)\cong L_\thick(\C)$, and therefore $L_\thick(\C)$ is an ideal
lattice.
\end{proof}

It would be interesting to have a direct proof (not involving a
classification) that the thick tensor ideals of $\coh X$ form an ideal
lattice. Note that Proposition~\ref{pr:abel} does not apply because the
tensor product $\otimes_{\OO_X}$ is exact only in trivial cases.

\begin{thm}\label{th:coh}
Let $X$ be a noetherian scheme and consider the abelian tensor
category $\coh X$ of coherent $\OO_X$-modules.  The pair $(X,\supp_X)$
is a classifying support datum on $\coh X$ and there is an induced
isomorphism
$$(X,\OO_X)\xto{\sim}(\Spec \coh X,\OO_{\coh X})$$ of ringed spaces.
\end{thm}
\begin{proof}
Let $\C=\coh X$.  It follows from well-known properties of the support
$\supp_X(x)$ that $(X,\supp_X)$ is a support datum on $\C$. Thus we
obtain a continuous map $f\colon X\to\Spec \C$ satisfying
$\supp_X(x)=f^{-1}(\supp(x))$ for each $x\in\C$, by
Theorem~\ref{th:final}. The classification of thick subcategories
of $\C$ from Proposition~\ref{pr:cohX} shows that the support datum
$(X,\supp_X)$ is classifying. Here we use in addition that the
underlying space of $X$ is spectral. It follows from
Corollary~\ref{co:class} that $f$ is a homeomorphism.

It remains to construct an isomorphism
$f^\sharp\colon\OO_\C\to f_*\OO_X$. Observe that for each open
$U\subseteq \Spec \C$, the restriction $\coh X\to\coh f^{-1}U$ induces an
equivalence $\C/\C_U\xto{\sim}\coh f^{-1}U$; see \cite[Prop.~VI.2]{G}.
Thus we obtain for $e=\OO_X$ an isomorphism
$$\OO_\C(U)=\End_{\C/\C_U}(e)\xto{\sim}\OO_X(f^{-1}U)$$ which yields
the isomorphism $f^\sharp\colon\OO_\C\xto{\sim}f_*\OO_X$.
\end{proof}

\subsection{Perfect complexes on a scheme}

We consider a quasi-compact and quasi-sep\-a\-rated scheme $X$ and its
triangulated tensor category $\bfD^\per (X)$ of perfect complexes with
tensor product $\otimes_{\OO_X}^\bfL$; see \cite[Sec.~2]{TT} for a
concise discussion of these concepts. For instance, every noetherian
scheme is quasi-compact and quasi-separated. Let us recall Thomason's
classification of thick tensor ideals.  Given $x\in\bfD^\per(X)$, we
write
$$\supp_X(x)=\{P\in X\mid x_P\neq 0\}.$$

\begin{prop}\label{pr:perX}
Let $X$ be a quasi-compact and quasi-separated scheme. The assignments
$$\bfD^\per(X)\supseteq
\D\mapsto\bigcup_{x\in\D}\supp_X(x)\quad\text{and}\quad X\supseteq
Y\mapsto \{x\in\bfD^\per(X)\mid\supp_X(x)\subseteq Y\}$$ induce
bijections between
\begin{enumerate} 
\item the set of all thick tensor ideals of $\bfD^\per(X)$, and
\item the set of all subsets $Y\subseteq X$ of the form
$Y=\bigcup_{i\in \Omega} Y_i$ with quasi-compact open complement
$X\setminus Y_i$ for all $i\in \Omega$.
\end{enumerate}
\end{prop}
\begin{proof}
See \cite[Thm.~4.1]{T}.
\end{proof}

We observe that the thick tensor ideals of $\bfD^\per(X)$ form an
ideal lattice by Proposition~\ref{pr:tria}. The following result shows
that a quasi-compact and quasi-separated scheme can be reconstructed
from the triangulated tensor category of perfect complexes; it is
a slight generalization of \cite[Thm.~6.3]{B} which assumes the scheme
to be topologically noetherian.

\begin{thm}
Let $X$ be a quasi-compact and quasi-separated scheme and consider the
triangulated tensor category $\bfD^\per(X)$ of perfect complexes on
$X$. The pair $(X,\supp_X)$ is a classifying support datum on
$\bfD^\per(X)$ and there is an induced isomorphism
$$(X,\OO_X)\xto{\sim}(\Spec \bfD^\per(X),\OO_{\bfD^\per(X)})$$ of
ringed spaces.
\end{thm}
\begin{proof}
The proof is essentially the same as that of Theorem~\ref{th:coh},
with $\C=\coh X$ replaced by $\C=\bfD^\per(X)$.  Note that the
assumption on $X$ implies that the underlying space is spectral. We
use the classification of thick tensor ideals from
Proposition~\ref{pr:perX}. For the equivalence
$\C/\C_U\xto{\sim}\bfD^\per(f^{-1}U)$, up to direct factors, when
$U\subseteq \Spec \C$ is quasi-compact open, we refer to
\cite[Sec.~5]{TT}.
\end{proof}

\section{A projective scheme}

There are interesting examples of triangulated tensor categories $\C$
where $(\Spec \C,\OO_\C)$ is actually a projective scheme. Here, we
present a general criterion which explains those examples. We fix a
triangulated tensor category $\C$ with tensor identity $e$. In
addition we assume that the tensor product is exact in each variable
and that the tensor category is suspended in the sense of \cite{S}.
Let us start with some preparation.  For $x,y\in\C$, we write
$$\Hom^*_\C(x,y)=\coprod_{n\in\bbZ}\Hom_\C(x,\Si^n y)$$ where $\Si$
denotes the suspension of $\C$.  The graded endomorphism ring
$\End^*_\C(x)$ acts on $\Hom^*_\C(x,y)$ from the right and
$\End^*_\C(y)$ acts from the left. We use the graded ring homomorphism
$$\p_x\colon\End^*_\C(e)\lto\End^*_\C(x),\quad\a\mapsto\a\otimes
x=x\otimes\a.$$ Note that $\End^*_\C(e)$ acts on $\Hom^*_\C(x,y)$ from
the right via $\p_x$ and from the left via $\p_y$, with
$$\a\cdot\b=(-1)^{|\a||\b|}\b\cdot\a$$ for homogeneous elements $\a\in
\End^*_\C(e)$ and $\b\in\Hom^*_\C(x,y)$. This follows from arguments 
similar to those in \cite{S}.  In particular, $\End^*_\C(e)$ is graded
commutative.

\subsection{Cohomological localization}
We need a basic result about the localization of triangulated
categories.  Under appropriate assumptions, we show that first taking
cohomology and then localizing is the same as first localizing and
then taking cohomology.  For a homogeneous element $\s\colon \Si^{n}
x\to x$ in $\End^*_\C(x)$, we denote by $x/\s$ its cofiber in $\C$.

\begin{prop}\label{pr:loc}
Let $\C$ be a triangulated category and $\D\subseteq \C$ a full
triangulated subcategory. Let $c\in \C$ be an object and $\p\colon
H\to\End_\C^*(c)$ be a graded ring homomorphism such that $H$ is
graded commutative.  Fix a subset $S\subseteq H$ of homogeneous
elements and consider for each $x\in\C$ the following
commutative diagram of canonical homomorphisms in the category of
graded $H$-modules.
$$\xymatrix{
\Hom_\C^*(c,x)\ar[r]^-\mu\ar[d]_-\pi&S^{-1}\Hom_\C^*(c,x)\ar[d]^-{S^{-1}\pi}\\
\Hom_{\C/\D}^*(c,x)\ar[r]^-\nu&S^{-1}\Hom_{\C/\D}^*(c,x)}$$
\begin{enumerate}
\item If $\{c/\p(\s)\mid\s\in S\}\subseteq\D$, then $\nu$ is an
isomorphism.
\item If $\D\subseteq \{x\in\C\mid S^{-1}\Hom_\C^*(c,x)=0\}$, then
${S^{-1}\pi}$ is an isomorphism.
\end{enumerate}
\end{prop}
\begin{proof}
We assume that $H$ is graded commutative because we need in (2) that
localization of graded $H$-modules with respect to $S$ is an exact functor.  

(1) Assume $\{c/\p(\s)\mid\s\in S\}\subseteq\D$. Let
$Q\colon\C\to\C/\D$ denote the quotient functor. Then $H$ acts on
$\Hom_{\C/\D}^*(c,x)$ via $Q$, and each $\s\in S$ acts invertibly
since $Q\p(\s)$ is invertible. Thus the canonical map
$$\Hom_{\C/\D}^*(c,x)\to S^{-1}\Hom_{\C/\D}^*(c,x)$$ is
invertible.

(2) Assume $\D\subseteq \{x\in\C\mid S^{-1}\Hom_\C^*(c,x)=0\}$.
We embed $\C$ into the category $\Mod\C$ of additive functors
$\C^\op\to\Ab$ via the Yoneda functor $$\C\lto\Mod\C,\quad
x\mapsto\Hom_\C(-,x).$$ Note that every cohomological functor
$F\colon\C\to\A$ into an abelian Grothendieck category $\A$ extends
uniquely to an exact and coproduct preserving functor $\bar
F\colon\Mod\C\to\A$; see \cite[Lem.~2.2]{K}. Now take the composition
$$\C\xto{\Hom^*_\C(c,-)}\Mod H\xto{S^{-1}}\Mod H$$ which
annihilates $\D$ by our assumption. We obtain the following
commutative diagram
$$\xymatrix{
\C\ar[rr]^-{F=\Hom^*_\C(c,-)}\ar[d]_-{Q}&&\Mod H\ar[d]^-{S^{-1}}\\
\C/\D\ar[rr]^-G&&\Mod H}$$
which can be extended to the following
commutative diagram of exact and coproduct preserving functors.
$$\xymatrix{ \Mod\C\ar[r]^-{\bar F}\ar[d]_-{Q^*}&\Mod
H\ar[d]^-{S^{-1}}\\ \Mod\C/\D\ar[r]^-{\bar G}&\Mod H}$$ Note that
$$\bar
F(M)\cong\coprod_nM(\Si^{-n}c)\cong\coprod_n\Hom_\C(\Hom_\C(-,\Si^{-n}c),M)$$
for $M$ in $\Mod\C$. The first isomorphism is clear for each
representable functor $M=\Hom_\C(-,x)$. Then observe that every
$M\in\Mod\C$ is a colimit of representable functors and both functors
preserve colimits. The second isomorphism follows from Yoneda's lemma.
The functor $Q^*$ has a right adjoint
$$Q_*\colon \Mod\C/\D\lto\Mod\C,\quad M\mapsto M\comp Q,$$ and the
adjunction morphism $Q^* Q_*M\to M$ is an isomorphism for all
$M\in\Mod\C/\D$, since $Q$ is a quotient functor. Now consider for
$x\in\C$ the adjunction morphism
$$\eta_x\colon\Hom_\C(-,x)\to Q_*Q^*\Hom_\C(-,x).$$ First observe that
$Q^*\eta_x$ is an isomorphism.  On the other hand, $\bar F\eta_x$
equals $\pi$ up to an isomorphism, since
\begin{align*}
\bar F( Q_*Q^*\Hom_\C(-,x))
&\cong\coprod_n\Hom_\C(\Hom_\C(-,\Si^{-n}c), Q_*Q^*\Hom_\C(-,x))\\
&\cong\coprod_n\Hom_{\C/\D}(Q^*\Hom_\C(-,\Si^{-n}c),Q^*\Hom_\C(-,x))\\
&\cong\coprod_n\Hom_{\C/\D}(\Hom_{\C/\D}(-,\Si^{-n}c),\Hom_{\C/\D}(-,x))\\
&\cong\coprod_n\Hom_{\C/\D}(\Si^{-n}c,x)\\
&=\Hom^*_{\C/\D}(c,x).
\end{align*}
Thus $S^{-1}\pi\cong S^{-1}(\bar F\eta_x)= \bar G (Q^*\eta_x)$ is an
isomorphism, and this finishes the proof.
\end{proof}

We formulate an immediate consequence.

\begin{cor}
Let $\C$ be a triangulated category. Let $c\in \C$ be an object
such that its graded endomorphism ring $\End_\C^*(c)$ is graded
commutative and fix a homogeneous prime ideal $\fp$.  Let
$$\C_\fp=\{x\in\C\mid\Hom_\C^*(c,x)_\fp=0\}.$$ Then we have a natural
isomorphism $\Hom_{\C/\C_\fp}^*(c,x)\xto{\sim}\Hom_\C^*(c,x)_\fp$ for
all $x\in\C$.
\end{cor}

\subsection{Cohomological support}
We keep fixed a triangulated tensor category $\C$ with tensor identity
$e$ and suppose that $H=\End_\C^*(e)$ is concentrated in non-negative
degrees.  We define the {\em cohomological support} of an object
$x\in\C$ as
$$\supp_H(x)=\{\fp\in\Proj H\mid\End_\C^*(x)_\fp\neq 0\},$$ where $H$
acts on $\End_\C^*(x)$ via the canonical ring homomorphism
$H\to\End_\C^*(x)$ taking an element $\a$ to $\a\otimes x=x\otimes\a$.
It is useful to observe that for each $\fp\in\Proj H$
\begin{equation}\label{eq:suppx}
\End_\C^*(x)_\fp= 0\quad\iff\quad \Hom_\C^*(c,x)_\fp=0\text{\ \ for all\ \ }
c\in\C.
\end{equation}


\subsection{A projective scheme}

We provide a criterion for $(\Spec \C,\OO_\C)$ to be a projective
scheme. We assume that the tensor product on $\C$ is commutative or
that $\C$ is generated by a single object as a triangulated
category. Thus the thick tensor ideals of $\C$ form an ideal lattice,
by Proposition~\ref{pr:tria}.  The following elementary observation
will be useful.

\begin{lem}\label{le:sind}
Let $\C$ be an abelian or triangulated tensor category and
$(X,\s)$ be a support datum on $\C$.  Then
$\C_0=\{x\in\C\mid\s(x)=\emptyset\}$ is a thick tensor ideal and
$(X,\s)$ induces a support datum $(X,\s')$ on the quotient $\C/\C_0$
such that $\s'(x)=\s(x)$ for all $x\in\C$.
\end{lem}

The following is the main result.

\begin{thm}\label{th:proj}
Let $\C$ be a triangulated tensor category with tensor identity $e$
and suppose that $H=\End_\C^*(e)$ is concentrated in non-negative
degrees.  Define for each $x\in\C$
$$\supp_H(x)=\{\fp\in\Proj H\mid\End_\C^*(x)_\fp\neq 0\}$$ and
suppose that $$\supp_H(x\otimes y)=\supp_H(x)\cap\supp_H(y)\quad
\text{for all}\quad x,y\in\C.$$ Then $\C_0=\{x\in\C\mid
\supp_H(x)=\emptyset\}$ is a thick tensor ideal of $\C$ and $(\Proj
H,\supp_H)$ induces a support datum on the quotient $\bar\C=\C/\C_0$. We
obtain an induced morphism
$$f\colon (\Proj H,\OO_H)\to (\Spec \bar\C,\OO_{\bar\C})$$ of ringed spaces
which induces a ring isomorphism $\OO_{\bar\C,f(\fp)}\xto{\sim}\OO_{H,\fp}$
for all $\fp\in\Proj H$. In particular, $f$ is an isomorphism if and
only if the support datum $(\Proj H,\supp_H)$ on $\bar\C$ is classifying.
\end{thm}
\begin{proof}
The condition $\supp_H(x\otimes y)=\supp_H(x)\cap\supp_H(y)$ implies
that $(\Proj H,\supp_H)$ is a support datum on $\C$, and therefore
also on $\bar\C$ by Lemma~\ref{le:sind}.  Thus we obtain a continuous map
$f\colon\Proj H\to\Spec \bar\C$ satisfying
$$\supp_H(x)=f^{-1}(\supp(x))$$ for all $x\in\bar\C$, by
Theorem~\ref{th:final}.

We need to construct a morphism of sheaves $f^\sharp\colon\OO_{\bar\C}\to
f_*\OO_H$. First observe that for each $x\in\bar\C$, the ring $H$ acts on
$\Hom^*_{\bar\C}(e,x)$ via the quotient functor $\C\to\bar\C$. In particular,
$$\Hom_\C^*(e,x)_\fp\xto{\sim}\Hom_{\bar\C}^*(e,x)_\fp
\quad\text{for}\quad\fp\in\Proj H$$ by Proposition~\ref{pr:loc}.  Now fix an
open subset $U\subseteq \Spec \bar\C$ and consider the composition of
the functors
$$F\colon\bar\C\xto{\Hom^*_{\bar\C}(e,-)}\Mod
H\xto{(-)\,\tilde{}}\Qcoh\Proj H\xto{(-)|_{f^{-1}(U)}}\Qcoh
f^{-1}(U).$$ Here, we denote for any $H$-module $M$ by $\tilde M$ its
associated sheaf. Note that the stalk of $\tilde M$ at a homogeneous
prime $\fp$ equals the degree zero part $M_{(\fp)}$ of the localized
module $M_\fp$. We claim that $F$ annihilates $\bar\C_U$. In fact,
$x\in\bar\C_U$ implies $f^{-1}(\supp(x))\cap f^{-1}(U)=\emptyset$ and
therefore $\supp_H(x)\cap f^{-1}(U)=\emptyset$. Thus
$\Hom_{\bar\C}^*(e,x)_{(\fp)}=0$ for all $\fp\in f^{-1}(U)$ and
therefore $Fx=0$. It follows that $F$ factors through
$\bar\C/\bar\C_U$ and induces a map
$\End_{\bar\C/\bar\C_U}(e)\to\OO_H(f^{-1}(U))$ which extends to a map
$\OO_{\bar\C}(U)\to\OO_H(f^{-1}(U))$. This yields the morphism of
sheaves $f^\sharp\colon\OO_{\bar\C}\to f_*\OO_H$.

Now fix a point $\fp\in\Proj H$. Then $f^\sharp$ induces a map
$f^\sharp_\fp\colon\OO_{\bar\C,f(\fp)}\to\OO_{H,\fp}$. We have an
isomorphism $\OO_{\bar\C,f(\fp)}\cong\End_{\bar\C/f(\fp)}(e)$ by
Lemma~\ref{le:stalk}. Next observe that
$$f(\fp)=\{x\in\bar\C\mid\End^*_{\C}(x)_\fp=0\}.$$ We have
$$\{e/\s\in\bar\C\mid\s\in H\setminus\fp\}\subseteq
f(\fp)\subseteq\{x\in\bar\C\mid\Hom^*_{\bar\C}(e,x)_\fp=0\}.$$
This follows from (\ref{eq:suppx}), and we obtain a second
isomorphism
$$\End_{\bar\C/f(\fp)}(e)\cong\End^*_{\bar\C}(e)_{(\fp)}\cong
\End^*_{\C}(e)_{(\fp)}=\OO_{H,\fp}$$
from Proposition~\ref{pr:loc}. We conclude that $f^\sharp_\fp$ is an
isomorphism. It follows that $f$ is an isomorphism of ringed spaces if
and only if the map $\Proj H\to\Spec \bar\C$ is a homeomorphism.  This
last condition is satisfied if and only if the support datum $(\Proj
H,\supp_H)$ is classifying, by Corollary~\ref{co:class}.
\end{proof}

Note that Theorem~\ref{th:proj} gives a partial answer to Balmer's
question when $(\Spec \C,\OO_\C)$ is a scheme \cite[Rem.~6.4]{B}. The
result is best illustrated by the following example from
representation theory; see \cite[Thm.~7.3]{FP} for an alternative
discussion.

\begin{exm}
Let $k$ be a field and let $A=kG$ be the group algebra of a finite
group $G$ or more generally a finite group scheme. We consider the
category $\mod A$ of finite dimensional $A$-modules and its bounded
derived category $\bfD^b(\mod A)$. The tensor product $\otimes_k$
on $\mod A$ induces a tensor product on $\bfD^b(\mod A)$ which is
exact in each variable.  The trivial representation $k$ is the tensor
identity and its graded endomorphism ring equals the group cohomology
ring
$$H=H^*(G,k)=\Ext^*_A(k,k).$$ Note that for $x\in\bfD^b(\mod A)$,
we have $\supp_H(x)=\emptyset$ if and only if $x$ belongs to the thick
tensor ideal $\bfD^\per(A)$ of perfect complexes. The composite
$$\mod A\xto{\inc}\bfD^b(\mod A)\xto{\can}\bfD^b(\mod A)/\bfD^\per(A)$$
induces an equivalence
$\umod A\xto{\sim}\bfD^b(\mod A)/\bfD^\per(A)$, where $\umod A$
denotes the stable module category of $A$; see for instance
\cite[Thm~2.1]{R}. The thick tensor ideals of $\umod A$ have been
classified, in case $G$ is a finite group and $k$ is algebraically
closed by Benson, Carlson, and Rickard in \cite[Thm.~3.4]{BCR}, and
for a finite group scheme over an arbitrary field by Friedlander and
Pevtsova in \cite[Thm.~6.3]{FP}. The classification implies that
$(\Proj H,\supp_H)$ is a classifying support datum on $\umod A$, and
therefore we have an isomorphism
$$(\Proj H^*(G,k),\OO_{H^*(G,k)})\xto{\sim} (\Spec \umod
kG,\OO_{\umod kG})$$ of ringed spaces by Theorem~\ref{th:proj}.
\end{exm}

The following example shows a triangulated tensor category which
arises in modular representation theory. The tensor product is not
necessarily commutative. This category can be used to define support
varieties of representations of finite dimensional algebras,
generalizing the classical case of a group algebra; see \cite{SS}

\begin{exm}
Let $A$ be a finite dimensional algebra over a field $k$ and let
$A^e=A\otimes_kA^\op$ be its enveloping algebra.  We consider the
category $\mod A^e$ of finite dimensional $A^e$-modules and the full
subcategory $\B$ of $A^e$-modules which are projective when restricted
to $A$ or $A^\op$. Note that $\B$ carries an exact structure which is
induced from the natural exact structure of $\mod A^e$. The inclusion
$\B\to\mod A^e$ induces a fully faithful exact functor
$\bfD^b(\B)\to\bfD^b(\mod A^e)$. The tensor product $\otimes_A$ on
$\mod A^e$ is exact in each variable when restricted to $\B$ and
induces therefore an exact tensor product on $\bfD^b(\B)$. The
$A^e$-module $A$ viewed as a complex concentrated in degree zero is a
tensor identity of $\bfD^b(\B)$. The tensor product restricts to a
tensor product on the thick subcategory $\C\subseteq\bfD^b(\B)$ which
is generated by $A$. We have therefore a triangulated tensor category
$(\C,\otimes_A,A)$ and the lattice of thick tensor ideals is an ideal
lattice, by Proposition~\ref{pr:tria}. Note that the tensor product of
$\C$ is not necessarily commutative. The graded endomorphism
ring of the tensor identity
$$\End^*_\C(A)=\Ext^*_{A^e}(A,A)= \HH^*(A)$$ equals
the Hochschild cohomology ring of $A$.
\end{exm}

\section{Decompositions of ideals}

In this final section we sketch how the decomposition of objects and
ideals of an additive tensor category are reflected by the
decomposition of their supports. The prototypical result in this
direction is Carlson's theorem from modular representation theory
which says that the variety of an indecomposable module is connected
\cite{Ca}.

\subsection{Decompositions of ideals}

Let $L$ be an ideal lattice and write $0=\inf L$. A non-zero element
$a\in L$ is called {\em indecomposable} if $a=a_1\vee a_2$ implies
$a_1=0$ or $a_2=0$.

\begin{prop}\label{pr:decomp}
Let $L$ be an ideal lattice and suppose that the space
$\Spec^*  L$ is noetherian. Given a semi-prime $a\in L$, there exists a
unique decomposition $a=\bigvee_{i\in \Om} a_i$ such that
\begin{enumerate} 
\item $a_i$ is indecomposable and semi-prime for all $i\in \Om$, and
\item $a_i\wedge a_j=0$ for all $i\neq j$ in $\Om$.
\end{enumerate}
\end{prop}
\begin{proof}
  First observe that every open subset of $\Spec^*  L$ is quasi-compact
  since $\Spec^*  L$ is noetherian. Thus Proposition~\ref{pr:bijection*}
  provides a bijection $b\mapsto \supp(b)$ between all semi-primes in
  $L$ and all subsets of $\Spec^*  L$ which are unions of closed
  subsets.  Under this bijection, a decomposition $a=\bigvee_i a_i$
  satisfying (1) and (2) corresponds to a disjoint union
  $$\supp(a)=\bigcup_i \supp(a_i).$$ Now observe that the unions of
  closed subsets are closed under arbitrary intersections. Thus there
  exists a partition $\supp(a)=\bigcup_i Y_i$ into unions of closed
  subsets which admits no proper refinement; see Lemma~\ref{le:part}
  below. We obtain the decomposition $a=\bigvee_i a_i$ by taking for
  $a_i$ the semi-prime satisfying $\supp(a_i)=Y_i$.
\end{proof}

\begin{lem}\label{le:part}
Let $X$ be a set and $\Y$ be a family of subsets which is closed under
forming intersections. Then there exists for each $Y\in\Y$ a unique
partition $Y=\bigcup_{i\in\Om} Y_i$ into non-empty subsets from $\Y$ which
admits no proper refinement.  More precisely, for all $i$, a disjoint
union $Y_i=Y_{i1}\cup Y_{i2}$ with $Y_{i1},Y_{i2}\in \Y$ implies
$Y_{i1}=\emptyset$ or $Y_{i2}=\emptyset$.
\end{lem}
\begin{proof}
Let $(\bigcup_{i\in\Omega_s}Y_{si})_{s\in\Si}$ be the family of all
partitions of $Y$ with $Y_{si}\in\Y$ for all $s,i$. For each $x\in Y$,
let
$$Y_x=\bigcap_{\genfrac{}{}{0pt}{}{s\in\Si}{x\in Y_{si}}}Y_{si}.$$ Then
$Y=\bigcup_{x\in Y} Y_x$ is a partition which admits no proper
refinement.
\end{proof}

\begin{rem}
There are refinements of Proposition~\ref{pr:decomp} which do not
require the space $\Spec^* L$ to be noetherian. For instance, if $a\in
L$ is compact and the space $\supp (a)$ is noetherian (with the
induced topology), then we have a unique decomposition
$a=\bigvee_{i=1}^n a_i$ into indecomposables such that $a_i\wedge
a_j=0$ for all $i\neq j$.
\end{rem}

\subsection{Decompositions in tensor categories}

Let $\C$ be an abelian or triangulated tensor category. We consider
the lattice $L_\thick(\C)$ of thick tensor ideals of $\C$ and recall
the following definition from \cite{K1}. Given a thick tensor ideal
$\D$, a family $(\D_i)_{i\in\Om}$ of thick tensor ideals is a {\em
decomposition} of $\D$ if
\begin{enumerate}
\item the objects in $\D$ are the finite coproducts of objects from
the $\D_i$, and
\item $\D_i\cap\D_j=0$ for all $i\neq j$.
\end{enumerate}
Such a decomposition is denoted by $\D=\coprod_{i\in\Om}\D_i$, and we
call $\D$ {\em indecomposable} if $\D\neq 0$ and any decomposition
$\D=\D_1\amalg\D_2$ implies $\D_1=0$ or $\D_2=0$.

The decomposition $\D=\bigvee_i\D_i$ of a thick tensor ideal (as
discussed in Proposition~\ref{pr:decomp}) amounts to a
decomposition $\D=\coprod_i\D_i$, provided that
$$\D_1\wedge
\D_2=0\quad\Longrightarrow\quad\D_1\vee\D_2=\D_1\amalg\D_2$$ for every
pair of thick tensor ideals $\D_1,\D_2$. This property holds if $\C$
admits an appropriate internal Hom-functor, because then $\D_1\wedge
\D_2$ implies $\Hom_\C(\D_1,\D_2)=0$. We do not go into details, but
refer to the literature.  A treatment of decompositions in the stable
module category $\umod kG$ can be found in \cite{K1}, where $kG$
denotes the group algebra of a finite group $G$. For further
discussions, see the recent work of Balmer \cite{B2} and Chebolu
\cite{C}.

\end{document}